\documentclass[12pt]{article}
\usepackage{amssymb}
\usepackage{amsfonts}
\usepackage{times}
\usepackage{mathptmx}
\usepackage{amsmath}
\usepackage[usenames]{color}
\usepackage{mathrsfs}
\usepackage{amsfonts}
\usepackage{amssymb,amsmath}
\usepackage{CJK}
\usepackage{cite}
\usepackage{cases}
\usepackage{amsthm}
\def\cl{\centerline}

\def\vs{\vspace*}

\def\Z{\mathbb{Z}}

\def\C{\mathbb{C}}

\def\QED{\hfill$\Box$}
\def\ni{\noindent}

\numberwithin{equation}{section}
\newtheorem{theo}{Theorem}[section]
\newtheorem{defi}[theo]{Definition}

\newtheorem{lemm}[theo]{Lemma}

\begin{document}
\begin{center}
{\bf\large Super-biderivations of Lie superalgebras}
\footnote {Supported by the National Natural Science Foundation of China (No. 11431010, 11371278,  11161010).

$^{\,\dag}$Corresponding author: G.~Fan:~yzfanguangzhe@126.com.
}
\end{center}

\cl{Guangzhe Fan$^{\,*,\,\dag}$, Xiansheng Dai$^{\,*}$
}

\cl{\small $^{\,*}$Department of Mathematics, Tongji University, Shanghai 200092,
P. R. China}

\vs{8pt}

{\small\footnotesize
\parskip .005 truein
\baselineskip 3pt \lineskip 3pt
\noindent{{\bf Abstract:} In this paper we attempt to investigate the super-biderivations of Lie superalgebras. Furthermore, we prove that all super-biderivations on the centerless super-Virasoro algebras are inner super-biderivations. Finally, we study the linear super commuting maps on the centerless super-Virasoro algebras.
\vs{5pt}


}}
\parskip .001 truein\baselineskip 6pt \lineskip 6pt

\section{Introduction}
Throughout this paper, we denote by $\C,\, \Z$ the sets of complex numbers, integers respectively. We assume that all vector spaces are based on $\C$, unless otherwise stated.

As is well known, derivations and generalized derivations are very important subjects in the research of both algebras and their generalizations. In recent years, biderivations have aroused many scholars' great interests in \cite{B,B1,B2,B3,WY,WYC,ZFLW,C}. In \cite{B1,B2}, Bre$\check{s}$ar et al showed that all biderivations on commutative prime rings are inner bidrivations and they determined the biderivations of semiprime rings. The notation of biderivations of Lie algebras was introduced in \cite{WYC}. In addition, in \cite{WY} the authors proved that the biderivations on the Schr$\ddot{\rm o}$dinger-Virasoro algebra are inner biderivations. Furthermore, in \cite{C} the author obtained every biderivation of a simple generalized Witt algebra over a field of characteristic 0 is a inner biderivation. It may be useful and interesting for computing the biderivations of some important Lie (super)algebras.

Lie superalgebras as a generalization of Lie algebras came from supersymmetry in mathematical physics. The theory of Lie superalgebras plays prominent roles in modern mathematics and physics. In the present paper, we attempt to research the biderivations of Lie superalgebras. We hope that biderivations would help to promote the development of structure theories of Lie superalgebras. This is our motivation to present this paper.

Super-Virasoro algebras are closely related to the conformal field theory and the string theory. They play very important roles in many mathematics and physics.
For the last few years, the structures theories and representation theories of super-Virasoro algebras have been extensively undertaken by many authors in \cite{SZ,Y,Z}. The centerless super-Virasoro algebras $S(s):=SVir[s](s=0~or~\frac{1}{2})$ are Lie superalgebras whose even part $S_{\bar{0}}(s)$ is spanned by $\{L_{m}~|~m\in \Z\}$ and odd part $S_{\bar{1}}(s)$ is spanned by $\{G_{k}~|~k\in s+\Z\}$, equipped with the following Lie super brackets:
\begin{equation*}\label{100}
\aligned
&[L_{m},L_{n}]=(n-m)L_{m+n},\\
&[L_{m},G_{k}]=(k-\frac{m}{2})G_{m+k},\\
&[G_{k},G_{l}]=2L_{k+l},
\endaligned
\end{equation*}
for all $m,n\in \Z$, $k,l\in s+\Z$. 

This paper is organized as follows. In Section $2$, we review the basic notations which are used in this paper. In section $3$, we obtain some useful conclusions about the super-biderivations of Lie superalgebras. In Section $4$, we compute the super-biderivations on the centerless super-Virasoro algebras.
Furthermore, we have proved that all super-biderivations of the centerless super-Virasoro algebras are inner super-biderivations. Finally, we investigate the linear super commuting maps on the centerless super-Virasoro algebras.

\section{Preliminaries}
In this section, for the reader's convenience, we shall summarize some basic facts about Lie superalgebras used in this paper, see \cite{CW,K,S}. For convenience, the degree of $x$ or $\phi$ is denoted by $|x|$ or $|\phi|$.

\begin{defi}\label{2000}\rm
A Lie superalgebra is a superalgebra $L=L_{\bar{0}}\oplus L_{\bar{1}}$ with multiplication $[\cdot,\cdot]$ satisfying the following two axioms:
\begin{equation}\label{200}
\aligned
&skew-supersymmetry:~~~~[x,y]=-(-1)^{|x||y|}[y,x],\\
&super~~Jacobi~~identity:~~~~[x,[y,z]]=[[x,y],z]+(-1)^{|x||y|}[y,[x,z]],
\endaligned
\end{equation}
for all homogeneous $x,y,z\in L$.
\end{defi}

Obviously, we know that $L_{\bar{0}}$ is a Lie algebra. Thus, if $L_{\bar{1}}=0$, then $L$ is just a Lie algebra.

Next we shall recall some conclusions about biderivations of Lie algebras in \cite{WY,WYC,C}. There was a little problem about some conclusions in \cite{WY}. They ignored that biderivations of Lie algebras should be skew-symmetry when studying some properties of biderivations. We present some modified results of biderivations of Lie algebras in \cite{C}. For convenience, we assume that all biderivations on Lie (super)algebras are skew-(super)symmetry in this paper.

\begin{defi}\rm
Let $G$ be a Lie algebra, we call a bilinear map $\phi$:$G\times G\longrightarrow G$ a biderivation of $G$ if it satisfies the following conditions:
\begin{equation}\label{201}
\aligned
&skew-symmetry:~~~~\phi(x,y)=-\phi(y,x),
\endaligned
\end{equation}
\begin{equation}\label{202}
\aligned
&\phi([x,y],z)=[x,\phi(y,z)]+[\phi(x,z),y],
\endaligned
\end{equation}
\begin{equation}\label{203}
\aligned
&\phi(x,[y,z])=[\phi(x,y),z]+[y,\phi(x,z)],
\endaligned
\end{equation}
for all $x,y,z\in G$.
\end{defi}

\begin{defi}\label{210}\rm
Let $\lambda\in \C$, then the map $\phi_{\lambda}$:$G\times G\longrightarrow G$, sending $(x,y)$ to $\lambda[x,y]$, is a biderivation of $G$. All biderivations of this kind are called inner biderivations of $G$.
\end{defi}

\begin{lemm}\label{211}\rm
Let $\phi$ be a biderivation on $G$, then $[\phi(x,y),[u,v]]=[[x,y],\phi(u,v)]$ for all $x,y,u,v\in G$. In particular, we have $[\phi(x,y),[x,y]]=0$ for all $x,y\in G$.
\end{lemm}

\begin{lemm}\label{212}\rm
Let $\phi$ be a biderivation on $G$. If $[x,y]=0$, then $\phi(x,y)\in Z([G,G])$, where $Z([G,G])$ is the center of $[G,G]$.
\end{lemm}

\section{Super-biderivations}
In \cite{XW}, the authors mainly discussed some properties of super-biderivations on Heisenberg superalgebras. However, only when the biderivations should be skew-supersymmetry, their conclusions would be right. Hence, we modify some results about super-biderivations of Lie superalgebras.
In this section, we are concerned with the super-biderivations of Lie superalgebras. We suppose that $L$ is a Lie superalgebra, unless otherwise stated.

\begin{defi}\label{3000}\rm
We call a bilinear map $\phi$:$L\times L\longrightarrow L$ a super-biderivation of $L$ if 
it satisfies the following three equations:
\begin{equation}\label{301}
\aligned
&skew-supersymmetry:~~~~\phi(x,y)=-(-1)^{|x||y|}\phi(y,x),
\endaligned
\end{equation}
\begin{equation}\label{302}
\aligned
&\phi([x,y],z)=(-1)^{|\phi||x|}[x,\phi(y,z)]+(-1)^{|y||z|}[\phi(x,z),y],
\endaligned
\end{equation}
\begin{equation}\label{303}
\aligned
&\phi(x,[y,z])=[\phi(x,y),z]+(-1)^{(|\phi|+|x|)|y|}[y,\phi(x,z)],
\endaligned
\end{equation}
for all homogeneous $x,y,z\in L$.
\end{defi}

\begin{defi}\label{304}\rm
A super-biderivation $\phi$ of homogenous $\gamma\in\Z_2$ of $L$  is a
super-biderivation such that $\phi(L_{\alpha},L_{\beta})\subseteq L_{\alpha+\beta+\gamma}$ for any $\alpha,\beta\in\Z_2$.
Denote by ${\rm BDer}_{\gamma}(L)$ the set of all super-biderivations of homogenous $\gamma$ of $L$.
Obviously, ${\rm BDer\,}({L})={\rm BDer\,}_{\bar{0}}({L})\oplus{\rm BDer\,}_{\bar{1}}({L})$.
\end{defi}

\begin{lemm}\label{3200}\rm
If the map $\phi_{\lambda}$:$L\times L\longrightarrow L$, defined by $\phi_{\lambda}(x,y)=\lambda[x,y]$ for all homogenous
$x,y\in L$, where $\lambda\in \C$, then $\phi_{\lambda}$ is a super-biderivation of $L$. We call this class super-biderivations by inner super-biderivations.
\end{lemm}
\ni\ni{\it Proof.}\ \ Firstly, according to $\phi_{\lambda}(x,y)=\lambda[x,y]$, we deduce that $\phi_{\lambda}$ is an even super-biderivation, i.e., $|\phi_{\lambda}|=\bar{0}$.
Due to the skew-symmetry of Lie superalgebras, it can be readily shown that $$\phi_{\lambda}(x,y)=-(-1)^{|x||y|}\phi_{\lambda}(y,x)$$
for any homogenous $x,y\in L$.

Now we can easily obtain the following three equalities:
\begin{equation}\label{305}
\aligned
&\phi_{\lambda}([x,y],z)=\lambda[[x,y],z],
\endaligned
\end{equation}

\begin{equation}\label{306}
\aligned
&(-1)^{|\phi_{\lambda}||x|}[x,\phi_{\lambda}(y,z)]=\lambda[x,[y,z]],
\endaligned
\end{equation}

\begin{equation}\label{307}
\aligned
&(-1)^{|y||z|}[\phi_{\lambda}(x,z),y]=(-1)^{|y||z|}\lambda[[x,z],y].
\endaligned
\end{equation}

Using the super-Jacobi identity $[[x,y],z]=[x,[y,z]]+(-1)^{|y||z|}[[x,z],y]$,
thus we get $$\phi_{\lambda}([x,y],z)=(-1)^{|\phi_{\lambda}||x|}[x,\phi_{\lambda}(y,z)]+(-1)^{|y||z|}[\phi_{\lambda}(x,z),y]$$ for any homogeneous $x,y,z\in L$.

Similarly, we deduce that $\phi_{\lambda}(x,[y,z])=[\phi_{\lambda}(x,y),z]+(-1)^{(|\phi_{\lambda}|+|x|)|y|}[y,\phi_{\lambda}(x,z)]$ for any homogeneous $x,y,z\in L$.

Therefore, $\phi_{\lambda}$ is a super-biderivation.\QED

\begin{lemm}\label{3300}\rm
Let $\phi$ be a super-biderivation on $L$, then $[\phi(x,y),[u,v]]=(-1)^{|\phi|(|x|+|y|)}[[x,y],\phi(u,v)]$ for any homogeneous $x,y,u,v\in L$.
\end{lemm}
\ni\ni{\it Proof.}\ \ Firstly, according to the definitation of super-biderivations, we compute $\phi([x,u],[y,v])$ by different ways.

On the one hand, one has
\begin{eqnarray*}\label{309}
\phi([x,u],[y,v])&=&(-1)^{|\phi||x|}[x,\phi(u,[y,v])]+(-1)^{|u|(|y|+|v|)}[\phi(x,[y,v]),u]\\
&=&(-1)^{|\phi||x|}([x,[\phi(u,y),v]]+(-1)^{(|\phi|+|u|)|y|}[x,[y,\phi(u,v)]])\\
&&+(-1)^{|u|(|y|+|v|)}([[\phi(x,y),v],u]+(-1)^{(|\phi|+|x|)|y|}[[y,\phi(x,v)],u]).
\end{eqnarray*}

On the other hand, it follows that
\begin{eqnarray*}\label{310}
\phi([x,u],[y,v])&=&[\phi([x,u],y),v]+(-1)^{|y|(|\phi|+|x|+|u|)}[y,\phi([x,u],v)]\\
&=&(-1)^{|\phi||x|}[[x,\phi(u,y)],v]+(-1)^{|u||y|}[[\phi(x,y),u],v]\\
&&+(-1)^{|y|(|\phi|+|x|+|u|)}((-1)^{|\phi||x|}[y,[x,\phi(u,v)]]+(-1)^{|u||v|}[y,[\phi(x,v),u]]).
\end{eqnarray*}

Comparing two sides of the above two equations, and using the super Jacobi identity of Lie superalgebras, we obtain
\begin{equation*}
[\phi(x,y),[u,v]]-(-1)^{|\phi|(|x|+|y|)}[[x,y],\phi(u,v)]
\end{equation*}
\begin{equation}\label{311}
=(-1)^{|u||v|+|v||y|+|y||u|}([\phi(x,v),[u,y]]-(-1)^{|\phi|(|x|+|v|)}[[x,v],\phi(u,y)]).
\end{equation}

Now let $\Phi(x,y;u,v)=[\phi(x,y),[u,v]]-(-1)^{|\phi|(|x|+|y|)}[[x,y],\phi(u,v)]$.

According to the equality \eqref{311}, it suggests that $$\Phi(x,y;u,v)=(-1)^{|u||v|+|v||y|+|y||u|}\Phi(x,v;u,y).$$

For one thing, we get
\begin{eqnarray*}\label{312}
\Phi(x,y;u,v)&=&-(-1)^{|u||v|}\Phi(x,y;v,u)\\
&=&-(-1)^{|u||v|}(-1)^{|u||v|+|v||y|+|y||u|}\Phi(x,u;v,y)\\
&=&(-1)^{|y||u|}\Phi(x,u;y,v).
\end{eqnarray*}

For another, we also have
\begin{eqnarray*}\label{313}
\Phi(x,y;u,v)&=&(-1)^{|u||v|+|v||y|+|y||u|}\Phi(x,v;u,y)\\
&=&-(-1)^{|u||v|+|v||y|+|y||u|}(-1)^{|y||u|}\Phi(x,v;y,u)\\
&=&-(-1)^{|y||u|}(-1)^{2(|u||v|+|v||y|+|y||u|)}\Phi(x,u;y,v)\\
&=&-(-1)^{|y||u|}\Phi(x,u;y,v).
\end{eqnarray*}

Hence, it follows at once that
\begin{eqnarray}\label{315}
\Phi(x,y;u,v)=-\Phi(x,y;u,v).
\end{eqnarray}

Due to the field $\C$, this forces that
$$\Phi(x,y;u,v)=0,$$ i.e.,
\begin{eqnarray}\label{316}
[\phi(x,y),[u,v]]=(-1)^{|\phi|(|x|+|y|)}[[x,y],\phi(u,v)].
\end{eqnarray}

Hence, we have completed the proof.\QED

\begin{lemm}\label{3400}\rm
If $|x|+|y|=\bar{0}$, then $[\phi(x,y),[x,y]]=0$ for any homogenous $x,y\in L$.
\end{lemm}
\ni\ni{\it Proof.}\ \
Letting $u=x,v=y$ in \eqref{316}, then we get
\begin{eqnarray*}\label{317}
[\phi(x,y),[x,y]]&=&(-1)^{|\phi|(|x|+|y|)}[[x,y],\phi(x,y)]\\
&=&-(-1)^{2|\phi|(|x|+|y|)+2|x||y|+|x|^{2}+|y|^{2}}[\phi(x,y),[x,y]]\\
&=&-(-1)^{|x|^{2}+|y|^{2}}[\phi(x,y),[x,y]]\\
&=&-(-1)^{|x|+|y|}[\phi(x,y),[x,y]].
\end{eqnarray*}

In view of $|x|+|y|=\bar{0}$, then we have
$$[\phi(x,y),[x,y]]=-[\phi(x,y),[x,y]].$$

Therefore, it shows that
$$[\phi(x,y),[x,y]]=0.$$\QED

\begin{lemm}\label{3500}\rm
Let $\phi$ be a super-biderivation on $L$. If $[x,y]=0$, then $\phi(x,y)\in Z([L,L])$, where $Z([L,L])$ is the center of $[L,L]$.
\end{lemm}
\ni\ni{\it Proof.}\ \ If $[x,y]=0$, then we obtain
\begin{eqnarray*}\label{317}
[\phi(x,y),[u,v]]&=&(-1)^{|\phi|(|x|+|y|)}[[x,y],\phi(u,v)]\\
&=&0,
\end{eqnarray*}
for any $u,v\in L$.

Thus, we get that $\phi(x,y)$ commutes with $[L,L]$ i.e., $\phi(x,y)\in Z([L,L])$.\QED

\section{Super-biderivations of the centerless super-Virasoro algebras}
In this section, we would like to compute super-biderivations of the center super-Virasoro algebras.
\begin{lemm}\label{4000}\rm
Every super-biderivation on the centerless super-Virasoro algebra $S(0)$ is an inner super-biderivation.
\end{lemm}
\ni\ni{\it Proof.}\ \ Let $\phi$ be a super-biderivation on the centerless super-Virasoro algebra $S(0)$.

In view of Lemma~3.6, we get $\phi(L_{0},L_{0})\in Z([S(0),S(0)])$. Thus, $\phi(L_{0},L_{0})=0$.
Assume that $\phi(L_{0},L_{n})=\Sigma_{m\in\Z}(a^{n}_{m}L_{m}+b^{n}_{m}G_{m})$, where $a^{n}_{m},b^{n}_{m}\in\C$ for any $m\in\Z$.

According to $L_{m}\in S_{\bar{0}}(0)$, it follows that $|L_{m}|+|L_{n}|=\bar{0}$ for any $m,n\in\Z$. Based on Lemma~3.5, then we have $$[[L_{0},L_{n}],\phi(L_{0},L_{n})]=0.$$

Furthermore, we have seen that $$n[L_{n},\Sigma_{m\in\Z}(a^{n}_{m}L_{m}+b^{n}_{m}G_{m})]=0.$$

Moreover, we get the following two equalities:
$$n a^{n}_{m}(m-n)=0,$$
$$n b^{n}_{m}(m-\frac{n}{2})=0.$$

This implies that if $m\neq n$, then $a^{n}_{m}=0$.
We also have if $m\neq \frac{n}{2}$, then $b^{n}_{m}=0$.

Therefore,
\begin{equation}\label{403}
\phi(L_{0},L_{n})=\left\{\begin{array}{lll}
a^{n}_{n}L_{n}+b^{\frac{n}{2}}_{\frac{n}{2}}G_{\frac{n}{2}}&\mbox{if \ }n~is~even,\\[4pt]
a^{n}_{n}L_{n}&\mbox{if \ }n~is~odd.
\end{array}\right.
\end{equation}
Obviously, we know that $a^{0}_{0}=b^{0}_{0}=0$.

From Lemma~3.4, one has
\begin{equation}\label{404}
[\phi(L_{0},L_{n}),[L_{0},L_{1}]]=(-1)^{|\phi|(|L_{0}|+|L_{n}|)}[[L_{0},L_{n}],\phi(L_{0},L_{1})].
\end{equation}

This implies that
\begin{equation}\label{405}
a^{n}_{n}(n-1)=n a^{1}_{1}(n-1),
\end{equation}
for any $n\in\Z$. Hence, we deduce that $a^{n}_{n}=n a^{1}_{1}$.

In addition, we also obtain that when $n$ is even, then $\frac{n-1}{2}b^{\frac{n}{2}}_{\frac{n}{2}}=0$. Thus, we get $b^{\frac{n}{2}}_{\frac{n}{2}}=0$ for any $n\in 2\Z$.

Letting $\lambda=a^{1}_{1}$, we can verify, for arbitrary $n\in\Z$, that
\begin{equation}\label{406}
\phi(L_{0},L_{n})=\lambda [L_{0},L_{n}].
\end{equation}

In the following we would like to determine $\phi(L_{0},G_{n})$.

Due to $[L_{0},G_{0}]=0$, then we have $\phi(L_{0},G_{0})=0$ by Lemma~3.6.
Assume that $\phi(L_{0},G_{n})=\Sigma_{m\in\Z}(c^{n}_{m}L_{m}+d^{n}_{m}G_{m})$, where $c^{n}_{m},d^{n}_{m}\in\C$ for any $m\in\Z$.

Set $x=L_{0},y=G_{n},z=G_{0}$ in \eqref{303}, then we conclude
\begin{equation}\label{407}
\aligned
&\phi(L_{0},[G_{n},G_{0}])=[\phi(L_{0},G_{n}),G_{0}]+(-1)^{(|\phi|+|L_{0}|)|G_{n}|}[y,\phi(L_{0},G_{0})],
\endaligned
\end{equation}

It follows at once that
\begin{equation}\label{408}
\aligned
&2n\lambda L_{n}=[\Sigma_{m\in\Z}(c^{n}_{m}L_{m}+d^{n}_{m}G_{m}),G_{0}],
\endaligned
\end{equation}
Therefore, we obtain $d^{n}_{n}=n\lambda$, $d^{m}_{n}=0$ if $m\neq n$. We also deduce that
if $n\neq 0$, then $c^{n}_{m}=0$.
Hence, we have $\phi(L_{0},G_{n})=c^{n}_{0}L_{0}+n\lambda G_{n}$.

Set $x=L_{0},y=G_{1},z=G_{1}$ in \eqref{303}, then we have
\begin{equation}\label{409}
\aligned
&\phi(L_{0},[G_{1},G_{1}])=[\phi(L_{0},G_{1}),G_{1}]+(-1)^{(|\phi|+|L_{0}|)|G_{1}|}[G_{1},\phi(L_{0},G_{1})],
\endaligned
\end{equation}
This shows that$$2\lambda L_{2}=c^{1}_{0}G_{1}+(-1)^{|\phi|}(2\lambda L_{2}-c^{1}_{0}G_{1})$$

Finally, we get $|\phi|=\bar{0}$. Thus, $\phi(L_{0},G_{m})\subseteq S_{\bar{1}}(0)$.

Therefore, we have $\phi(L_{0},G_{n})=n\lambda G_{n}$. This means that we have proved the following equation
\begin{equation}\label{410}
\phi(L_{0},G_{n})=\lambda [L_{0},G_{n}],
\end{equation}
for any $n\in\Z$.

From \eqref{406} and \eqref{410}, we can easily obtain the following equation
\begin{equation}\label{411}
\phi(L_{0},z)=\lambda [L_{0},z],
\end{equation}
for any $z\in S(0)$.

Due to $|\phi|=\bar{0}$, and according to Lemma~3.3, we obtain
$$[\phi(x,y),[L_{0},z]]=[[x,y],\phi(L_{0},z)].$$

Thanks to the equality \eqref{411}, we get
$[\phi(x,y)-\lambda [x,y],[L_{0},z]]=0$.
According to the arbitrary of $z$, then
$\phi(x,y)-\lambda [x,y]\in Z(S(0))=0$.

Thus, $\phi(x,y)=\lambda [x,y]$.\QED

Using the similar method of proving Lemma~4.1, we can obtain the following lemma.
\begin{lemm}\label{4100}\rm
Every super-biderivation on the centerless super-Virasoro algebra $S(\frac{1}{2})$ is an inner super-biderivation.
\end{lemm}

Now we get a main result of this section.
\begin{theo}\label{4300}\rm
Every super-biderivation on the centerless super-Virasoro algebras $S(s)(s=0~or~\frac{1}{2})$ is an inner super-biderivation.
\end{theo}

Now we introduce some notations about super commuting maps. Let $L$ be an arbitrary Lie superalgebra. A linear map $f$ is called a linear super commuting map of $L$ if $[f(x),x]$ for all $x\in L$.

Furthermore, if $f$ is a linear super commuting map on $L$, then
\begin{eqnarray*}\label{416}
[f(x+y),x+y]&=&[f(x),x]+[f(x),y]+[f(y),x]+[f(y),y]\\
&=&[f(x),y]+[f(y),x]\\
&=&0.
\end{eqnarray*}
Thus, we have
\begin{equation}\label{420}
[f(x),y]=(-1)^{|f||x|}[x,f(y)],
\end{equation}
for any homogenous $x,y\in L$.

\begin{theo}\label{4400}\rm
Every linear super commuting map $f$ on the centerless super-Virasoro algebras $S(s)(s=0~or~\frac{1}{2})$ has the following form $f(x)=\lambda x$ for any $x\in S(s)$, where $\lambda\in\C$.
\end{theo}
\ni\ni{\it Proof.}\ \ Now we want to construct a linear map $\Psi:S(s)\times S(s)\longrightarrow S(s)$ such that $\Psi$ is a super-biderivation of $S(s)$.
Define $\Psi(x,y)=[f(x),y]$, it is easy to see that $|\Psi|=|f|$.

In view of super Jacobi identity, it suffices to prove that
\begin{equation}\label{421}
\aligned
&\Psi(x,[y,z])=[\Psi(x,y),z]+(-1)^{(|\Psi|+|x|)|y|}[y,\Psi(x,z)],
\endaligned
\end{equation}

According to the equality \eqref{420}, we also have
\begin{equation}\label{422}
\aligned
&\Psi([x,y],z)=(-1)^{|\Psi||x|}[x,\Psi(y,z)]+(-1)^{|y||z|}[\Psi(x,z),y],
\endaligned
\end{equation}

Thus, $\Psi$ is a super-biderivation on $S(s)$.

According to Theorem~4.3, we obtain $\phi(x,y)=\lambda [x,y]$. This implies that $[f(x),y]=\lambda [x,y]$, i.e., $[f(x)-\lambda x,y]=0$ for arbitrary $y\in S(s)$.

Hence, $f(x)=\lambda x$ for all $x\in S(s)$.

Conversely, if $f(x)=\lambda x$ for all $x\in S(s)$, then $f$ is a linear super commuting map on $S(s)$.
\QED

\end{document}